\magnification=\magstep1
\input amstex
\documentstyle{amsppt}
\catcode`\@=11 \loadmathfont{rsfs}
\def\mycal{\mathfont@\rsfs}
\csname rsfs \endcsname \catcode`\@=\active

\vsize=7.5in

\topmatter
\title 
%On the number of disjoint singular MASAs  \\ in an ultraproduct II$_1$ factor 
On the singular abelian rank \\ of ultraproduct II$_1$ factors 
\\
\\
$\text{\it Dedicated to Jacques Dixmier on his 100th birthday}$ 
  \endtitle
 
\author Patrick Hiatt and Sorin Popa$^*$  \endauthor

\rightheadtext{Singular rank of ultraproduct factors}

\affil     {\it  University of California, Los Angeles} \endaffil

\address Math.Dept., UCLA, Los Angeles, CA 90095-1555, USA \endaddress
\email  popa\@math.ucla.edu, hiatt\@math.ucla.edu \endemail

\thanks $^*$ Supported in part by NSF Grant DMS-1955812  and the Takesaki Endowed Chair at UCLA \endthanks

\abstract  We prove that, under the continuum hypothesis $\frak c=\aleph_1$,  any ultraproduct 
II$_1$ factor $M= \prod_{\omega} M_n$ of separable finite factors $M_n$ 
contains more than $\frak c$ many mutually disjoint singular MASAs, in other words the {\it singular abelian rank of} $M$, $\text{\rm r}(M)$, is larger than $ \frak c$. Moreover, if the strong continuum hypothesis  $2^{\frak c}=\aleph_2$ is assumed, then ${\text{\rm r}}(M) = 2^{\frak c}$. More generally, these results hold true for any II$_1$ factor $M$ with unitary 
group of cardinality $\frak c$ that satisfies  
the bicommutant condition $(A_0'\cap M)'\cap M=M$, for all $A_0\subset M$ separable abelian. 
%Since MASAs of ultraproduct II$_1$ factors are purely non-separable, this shows that the {\it sans-rank of} $M$, ${\text{\rm r}}_{\text{\rm  ns}}(M)$, is equal to $2^{\frak c}$ as well. 
\endabstract

\endtopmatter

\document

\heading 0. Introduction \endheading

Following  Dixmier [D54], a maximal abelian $^*$-subalgebra (MASA) $A$ in a von Neumann algebra $M$ is called {\it 
singular} if the only unitary elements $u\in \Cal U(M)$ that normalize $A$ (i.e., $uAu^*=A$) are the unitaries in $A$. The existence of such MASAs 
in the hyperfinite II$_1$ factor $R$ in [D54] was a discovery that led  to many interesting developments and subsequent research (see e.g., [Pu60], 
[P81a], [P81b], [R91], [P16],  [HP17]). 

Most recently in this direction, the {\it singular abelian core} of a II$_1$ factor $M$  was defined in [BDIP23] as 
the (unique up to unitary conjugacy) maximal abelian $^*$-subalgebra $A\subset \Cal M=M\overline{\otimes} \Cal B(\ell^2K)$, with $|K|\geq 2^{|\Cal U(M)|}$, 
that's generated by finite projections of $\Cal M$, is singular in $1_A\Cal M1_A$ and is maximal in $\Cal M$ with respect to inclusion.  Also, the {\it singular abelian rank  
of} $M$ was defined as $\text{\rm r}(M):=Tr_{\Cal M}(1_A)$, viewed as a cardinality when infinite. Alternatively, r$(M)$ can be viewed as the 
``maximal number'' of disjoint singular MASAs (or pieces of it) in $M$. The {\it sans-core} and respectively 
{\it sans-rank} $\text{\rm r}_{\text{\rm ns}}(M)$ were defined in [BDIP23] in a similar way, by considering the maximal singular abelian purely non-separable core 
$A\subset \Cal M=M\overline{\otimes} \Cal B(\ell^2K)$ and respectively the semi-finite trace of its support in $\Cal M$. 

It was pointed out in [BDIP23] that by results in ([P13a], [P18]), for any separable II$_1$ factor $M$ one has $\text{\rm r}(M)=\frak c$ and that 
if $M$ is an ultraproduct II$_1$ factor, $M=\Pi_\omega M_n$, associated to a sequence  $M_n$ of separable II$_1$ factors and a free ultrafilter 
$\omega$ on $\Bbb N$, then by simply considering ultraproducts of singular MASAs of $M_n$ 
one obtains $\text{\rm r}(M)=\text{\rm r}_{\text{\rm ns}}(M)\geq \frak c$. But a more exact calculation of the singular abelian  rank 
of such $M$ was left open. 

We prove in this paper that if we assume the continuum hypothesis (CH), $\frak c=2^{\aleph_0}=\aleph_1$, 
then for any II$_1$ factor of the form $M=\Pi_\omega M_n$, with $M_n$ separable tracial factors with dim$(M_n)\rightarrow \infty$,  
one has $\text{\rm r}(M)=\text{\rm r}_{\text{\rm ns}}(M)\geq 2^{\frak c}$,  
and that if we further assume the strong continuum hypothesis (SCH), $2^{\frak c}=\aleph_2$, then we actually 
have equalities, $\text{\rm r}(M)=\text{\rm r}_{\text{\rm ns}}(M)=2^{\frak c}$ (see Theorem 2.1). Note that in particular this shows that, 
under CH,  an ultraproduct II$_1$ factor has many more singular MASAs than the ones arising as ultraproducts of MASAs. 

To do this calculation, we in fact only use the property of an ultraproduct II$_1$ factor $M=\Pi_\omega M_n$  that any copy $A_0\subset M$ of the separable diffuse 
abelian von Neumann algebra $L^\infty[0,1]$ satisfies the bicommutant condition $(A_0'\cap M)'\cap M=A_0$. When viewed as an abstract property of a II$_1$ factor 
$M$, we call this {\it property} U$_0$. 

We prove that, somewhat surprisingly, a II$_1$ factor $M$ has property U$_0$ if and only if it has {\it property} U$_1$, 
requiring that any isomorphism between two copies of $L^\infty[0,1]$ inside $M$ is implemented by a unitary in $M$ (see Theorem 1.2), and 
call a II$_1$ factor satisfying any of these equivalent properties a U-{\it factor}. 

We also relate 
properties U$_0$, U$_1$ with  the 
weaker property that any two copies of $L^\infty[0,1]$ inside $M$ are unitary conjugate, already considered in ([P13a], [P20]), and which we label 
here U$_2$. This property for $M$ implies for instance that $M$ is prime and has 
no Cartan subalgebras and that any MASA in $M$ is purely non-separable (see Proposition 1.4). 
Thus, for such factors one always has r$_{\text{\rm ns}}(M)=\text{\rm r}(M)$. 

So with this terminology, our main result (Theorem 2.1) shows that if $M$ is a U-factor with unitary group $\Cal U(M)$ 
having cardinality $|\Cal U(M)|=\frak c$, 
then with the CH  assumption we have $\text{\rm r}(M)\geq 2^{\frak c}$, with equality when SCH is assumed.  

We mention that 
Gao, Kunnawalkam Elayavalli, Patchell and Tan have recently been able to construct 
examples of II$_1$ U-factors $M$ with $|\Cal U(M)|=\frak c$ but which cannot be decomposed as an  ultraproduct of separable  
finite factors ([GKPT24]). 

%were recently able to construct U-factors that are not ultraproduct II$_1$ factors ([GKPT24]). 

Throughout this paper we will systematically use notations,  terminology and basic results from [P13b] (for all things concerning ultraproduct 
II$_1$ factors) and [P16] (for intertwining of subalgebras and disjointness in II$_1$ factors, in particular for MASAs, especially singular ones). Our work here has been 
especially motivated by remarks and considerations in ([BDIP23], notably Sections 2.3, 2.4 and the remarks therein). We comment at length about this 
in Section 3 of this paper.  

We are very grateful to Adrian Ioana and Stefaan Vaes for many useful comments on a preliminary  draft of this paper.

\vskip.05in

{\it Acknowledgement}. This paper is dedicated to Jacques Dixmier, whose 
%one of the founding fathers of operator algebras, on the occasion of his 100th anniversary. His  
seminal monographs [D57], [D64] and many pioneering contributions to operator algebras  
played a crucial role in the development of this  area. 
%, a direction of research that has become increasingly relevant in today's mathematics. 
May this be just a small token of gratitude, on the occasion of his 100th anniversary, 
in the name of the several generations of mathematicians who benefitted from his work over more than seven decades.

\heading 1.   Some abstract properties of ultraproduct II$_1$ factors  \endheading

While any separable approximately finite dimensional (AFD) tracial von Neumann algebra $(B_0, \tau)$ can be embedded into any II$_1$ factor $M$ ([MvN43]), 
when $M$ is an ultraproduct II$_1$ factor, $M=\Pi_\omega M_n$, such an embedding $(B_0, \tau) \hookrightarrow M$ follows even unique 
up to unitary conjugacy in $M$.  Also, any separable AFD subalgebra $B_0\subset M$ satisfies the bicommutant condition $(B_0'\cap M)'\cap M = B_0$  
(see e.g. Theorem 2.1 in [P13b]). 

In particular, the uniqueness of the embedding and the bicommutant property  
hold true when $(B_0, \tau)$ is the separable diffuse abelian von Neumann algebra $(L^\infty [0,1], \int \cdot \ \text{\rm d}\lambda)$.  
In this section we will consider these  two properties as abstract properties of a II$_1$ factor $M$ and prove that 
they are in fact equivalent. We also discuss the apriori weaker condition that any two copies of $L^\infty [0,1]$ inside $M$ are unitary conjugate. 

\vskip.05in

\noindent
{\bf 1.1. Definition}. Given a II$_1$ factor  $M$, we consider the following three properties:

\vskip.05in
U$_0$ Any separable  abelian von Neumann subalgebra $A_0\subset M$ satisfies the bicommutant property $(A_0'\cap M)'\cap M=A_0$.

\vskip.05in
U$_1$ Any trace preserving isomorphism between two separable diffuse abelian von Neumann subalgebras of $M$ is implemented by a unitary element in $M$; 

\vskip.05in
U$_2$ Any two separable diffuse abelian von Neumann subalgebras of $M$ are unitary conjugate;

\vskip.05in

For each $i=0,1,2$, we say that $M$ has {\it stable property} U$_i$, if $M^t$ satisfies U$_i$ for any $t>0$. 

\proclaim{1.2. Theorem}  Conditions $\text{\rm U}_0, \text{\rm U}_1$ for a  $\text{\rm II}_1$ factor $M$ are equivalent and they are both stable 
properties, i.e, if $M$ satisfies property $\text{\rm U}_i$, 
for some $i=0,1$,  then $M^t$ satisfies  it  for any $t>0$.  
\endproclaim
\noindent
{\it Proof}. Let us first show that U$_1$ is stable. 
So assume $M$ satisfies U$_1$. We first show that  $N=\Bbb M_n(M)$ satisfies U$_1$ as well. Let $A_1, A_2\subset N$ be separable diffuse abelian von Neumann algebras 
and $\theta: A_1 \simeq A_2$ an isomorphism preserving the trace on $N$. Then $A_1$ contains a partition of $1$ with projections $\{p^1_j\}_{j=1}^n$ of trace equal $1/n$.  
Let $p^2_j=\theta(p^1_j)$. 
By conjugating with appropriate unitaries $u_1, u_2 \in N$ we may assume $p^i_j = e_{jj}$, $1\leq j \leq n$, $i=1,2$, where $\{e_{ij} \mid 1\leq i,j\leq n\}\subset \Bbb M_n(\Bbb C)$,  
are the matrix units. Denoting by $\theta_j$ the restriction of $\theta$ to $A_1e_{jj}\simeq A_2e_{jj}$ and viewing them both as subalgebras in $M\simeq e_{jj}Ne_{jj}$, 
by the U$_1$ property for $M$ it follows that $\theta_j$ is implemented by $u_j\in e_{jj}Ne_{jj}$. But then $u=\sum_j u_j\in \Cal U (N)$ implements $\theta:A_1\simeq A_2$. 

We now show that if $p\in \Cal P(M)$ then $pMp$ satisfies U$_1$. If $A_1, A_2\subset pMp$ are separable diffuse abelian von Neumann algebras 
and $\theta: A_1 \simeq A_2$ an isomorphism preserving the trace on $pMp$, then there exist separable diffuse abelian von Neumann subalgebras 
$\tilde{A}_i \subset M$ such that $p\in \tilde{A}_i$,  and $\tilde{A}_ip=A_i$, $i=1,2$, as well as  a trace preserving isomorphism $\tilde{\theta}:\tilde{A}_1 \simeq \tilde{A}_2$ 
whose restriction to $A_1$ is equal to $\theta$. If $u\in \Cal U(M)$ implements $\tilde{\theta}$, then $up\in \Cal U(pMp)$ implements $\theta$.  Thus, U$_1$ is stable.

Let us now prove that conditions U$_0$, U$_1$ are equivalent. 
Let $A_0\subset M$ be a separable diffuse abelian von Neumann algebra. Denote $B=A_0'\cap M$ and $Z=B'\cap M$. Note that $Z=\Cal Z(B)$. Indeed, because 
any element in $M$ that commutes 
with all elements in $B=A_0'\cap M$ must in particular commute with $A_0$ so $B'\cap M \subset B$, 
which is equivalent to  $B'\cap M=\Cal Z(B)$. 

Assume $M$ satisfies U$_1$. If $Z\neq A_0$, then there exists a projection $p\in Z$ with $b=E_{A_0}(p)\neq p$. There exists a projection 
$q\in A_0$ majorized by the support $s=s(b)$ of $b$ such that $cq \leq qb (1-c) q$ for some $c>0$. Let $B_0=A_0 \vee \{p\}$. Thus, by replacing $p$ 
by $qp$ we may assume $p$ itself satisfies $cs\leq b=E_{A_0}(p) \leq (1-c)s$. Denote $B_0=A_0s\vee \{p\} \subset Zp$.  
Note that the inclusion $L^\infty X \simeq A_0s \subset B_0 \simeq L^\infty Y$ 
is given by a surjective measure preserving map $\alpha:Y \rightarrow X$ with two-points fiber $\forall t\in X$. Consider then the trace preserving embedding of 
$(B_0, \tau_{B_0})$ into a tracial von Neumann algebra $Q\simeq A_0s\overline{\otimes}R$, endowed with the trace $\tau_{A_0s}\otimes \tau_R$,  
such that $A_0s$ identifies with the center $\Cal Z(Q)=A_0s \otimes 1\simeq L^\infty X$ and such that when we view $p$ as a measurable field 
$p_t, t\in X$, with $p_t\in \Cal P(R)$, we have $\tau_{R}(p_t)=b_t$,  where $(b_t)_t=b$. 

Since $Q$ with its trace can be embedded into any II$_1$ factor, we can view it as a von Neumann subalgebra of $sMs$ and then by using U$_1$ 
for $A_0s\subset sMs$ we may assume the center of $Q$ coincides with $A_0s$ and $B_0$ with $A_0s\vee \{p\}$. So $1\otimes R$ is in the 
commutant of $A_0s$, and hence of $A_0$. Since $p\in Z$, we should thus have $1\otimes R$ commute with $p$. But by averaging $p$ over the unitaries in $1\otimes R$ 
we get $b$, which is not equal to $p$, a contradiction. 

Thus, we must have $(A_0'\cap M)'\cap M=A_0$, showing that U$_0$ is satisfied. 

Conversely, assume $M$ satisfies the bicommutant condition U$_0$.  Let $A_1, A_2\subset M^{1/2}$ be separable diffuse abelian and $\theta:A_1 \simeq A_2$ 
be an isomorphism preserving the restrictions of the trace on $M^{1/2}$ to $A_1, A_2$. Let $A=\{a e_{11}+ \theta(a)e_{22} \mid a\in A_1\}$ which we view 
as a (separable abelian diffuse) von Neumann subalgebra of $M=\Bbb M_2(M^{1/2})$. Then $(A'\cap M)'\cap M=A$ implies in particular 
that the projections $e_{11}, e_{22} \in A'\cap M$ are equivalent in $A'\cap M$, via some partial isometry $v=ue_{12}$ where $u$ 
is a unitary in $e_{11}Me_{11}=M^{1/2}$. But this means $\theta(a)=uau^*$ for any $a\in A_1$. 

We have thus proved that if $M$ satisfies U$_0$ then $M^{1/2}$ satisfies U$_1$. Since we already showed that  U$_1$ is a stable property, this implies $M$ 
satisfies U$_1$. Thus, U$_0$, U$_1$ are equivalent, and since U$_1$ was shown to be stable, U$_0$ follows stable as well. 

\hfill $\square$

\vskip.05in 
\noindent
{\bf 1.3. Definition}. We say that a II$_1$ factor  $M$ is a U-{\it factor} 
 if it satisfies the equivalent conditions U$_0$, U$_1$. 
 
 \vskip.05in
 
 We already mentioned that ultraproduct II$_1$ factors $M=\Pi_\omega M_n$ satisfy the bicommutant property U$_0$ and 
 the unique (up to unitary conjugacy) embedding property U$_1$. They are the typical examples of U-factors. 
 
 Since property U$_1$ for a II$_1$ factor $M$ trivially implies the unitary conjugacy of any 
 two copies of $L^\infty[0,1]$ inside $M$,  i.e., condition U$_2$, any U-factor satisfies U$_2$ as well. Condition U$_2$ was already 
 considered as an abstract property of II$_1$ factors in (Proposition 2.3 of [P13a]), where it was noticed that the arguments 
 in (Section 7 of [P81]), showing that an ultraproduct II$_1$ factor $M$ has no Cartan subalgebras and all its MASAs are purely non-separable,   
 only use the fact that $M$ satisfies condition U$_2$.  
It was further noticed in ([P20]) that U$_2$ factors are prime and have the property that the commutant of any 
separable abelian $^*$-subalgebra is of type II$_1$. 

We restate all these results here, including their proofs from ([P81], [P13a], [P20]), for the reader's convenience.

 \vskip.05in

\proclaim{1.4. Proposition ([P81], [P13a], [P20])} Assume a   $\text{\rm II}_1$ factor $M$ satisfies property $\text{\rm U}_2$ 
$($for instance, if $M$ is a $\text{\rm U}$-factor$)$. 
Then $M$  automatically satisfies the following properties: 
\vskip.05in

$(a)$ Given any MASA $A$ in $M$,  there exists a diffuse abelian von Neumann subalgebra $B_0\subset M$ orthogonal to $A$. 

$(b)$ Any separable abelian von Neumann subalgebra $A_0\subset M$ has type $\text{\rm II}_1$ relative commutant $A_0'\cap M$. 

$(c)$ Any MASA in $M$ is purely non-separable. 

$(d)$ $M$ has no Cartan MASA. 

$(e)$ $M$ is prime. 
\endproclaim 
\noindent 
{\it Proof}.  $(a)$ Let $A\subset M$ be a MASA. Let $D\subset A$ be a separable diffuse von Neumann subalgebra. Since any two separable diffuse abelian subalgebras 
in $M$ are unitary conjugate and since $M$ contains copies of the hyperfinite II$_1$ factor (by [MvN43]),  we may assume $D$ is the 
Cartan subalgebra of such a subfactor $R \subset M$, 
represented as $D=D_2^{\otimes \infty}\subset M_{2 \times 2}(\Bbb C)^{\otimes \infty} =R$.  
Let $D_2^0\subset M_{2 \times 2}(\Bbb C)$ be a maximal abelian subalgebra of $M_{2 \times 2}(\Bbb C)$ that is perpendicular to $D_2$ 
and denote $D^0={D^0_2}^{\otimes \infty} \subset R$. Then $D \perp D^0$ and since both $D, D^0$ are MASAs in $R$, 
we have $E_{D'\cap M}(D^0)=E_{D'\cap R}(D^0)=E_{D}(D^0)=\Bbb C$, i.e. $D^0 \perp D'\cap M \supset A$, proving $1^\circ$. 

$(b)$ By [MvN43],  one has $R\overline{\otimes}R \simeq R$ and so $R\overline{\otimes}R$ embedds into $M$. If one takes any 
MASA $B_0\subset R\otimes 1\subset R\overline{\otimes}R \simeq R$, then $B_0'\cap M \supset 1\otimes R$, implying that 
$B_0'\cap M$ is type II$_1$. Since $A_0, B_0$ are unitary conjugate in $M$, $A_0'\cap M$ is II$_1$ as well.

$(c)$  Let $A$ be a MASA in $M$. If $Ap$ is separable for some projection $p\in M$, then by taking a smaller $p$ if necessary we may assume $\tau(p)=1/n$ 
for some integer $n\geq 1$. Let $v_1=p, v_2, ..., v_n\in M$ be partial isometries with $v_i^*v_i=p$, $\forall 1\leq i\leq n$, and $\sum_i v_iv_i^*=1$ and define 
$B=\sum_i v_i (Ap)v^*_i$. Then $B$ is a separable MASA in $M$. But then taking $B_0\subset B$ to be 
any diffuse proper von Neumann subalgebra of $B$,  
it cannot be unitary conjugate to $B$ because  $B_0$ is not a MASA while $B$ is, contradiction. 

$(d)$ Let $A\subset M$ be a MASA. By part $1^\circ$, there exist separable diffuse abelian subalgebras    
$D, D^0$ in $M$ such that $D\subset A$ and $D^0\perp A$. Let $u\in \Cal U(M)$ be so that $uDu^*=D^0$. Then $u$ is perpendicular 
to the normalizer of $A$ in $M$. Indeed, because for any $v\in \Cal N_M(A)$ and any partition $p_i\in D$ of mesh $\leq \varepsilon$, 
we have
$$
|\tau(uv)|^2=|\tau(\Sigma_i p_i uv p_i)|^2  
\leq \|\Sigma_i p_i u v p_i\|^2_2 = \Sigma_i \tau(u^*p_iuvp_iv^*) = \Sigma_i \tau(p_i)^2 
\leq \varepsilon. 
$$ 
Since $\varepsilon >0$ was arbitrary, $\tau(uv)=0$. Thus $u \perp \Cal N_M(A)''$. 

$(e)$ If $M=M_1 \overline{\otimes} M_2$ with $M_1, M_2$ of type II$_1$ then there exist separable diffuse abelian von Neumann subalgebras 
$A_i\subset M_i$. By hypothesis, there exists a unitary $u\in M$ such that $uA_1u^*= A_2\perp A_1$. 
From the argument in $4^\circ$, it follows that for any unitaries 
$v_1\in M_1$, $v_2\in M_2$ one has $\tau(uv_1v_2)=\tau(v_2uv_1)=0$. Taking span of $v_i$ 
and using that the $\| \ \|_2$ closure of the span of $1\otimes M_2 \cdot M_1\otimes 1$ is $M$, it follows that $\tau(uu^*)=0$, contradiction.  
\hfill $\square$

\proclaim{1.5. Corollary} If a $\text{\rm II}_1$ factor $M$ satisfies property $\text{\rm U}_2$ $($e.g., if $M$ is a $\text{\rm U}$-factor$)$, 
then $\text{\rm r}_{\text{\rm ns}}(M)=\text{\rm r}(M)$. 
\endproclaim
\noindent
{\it Proof}. By part $(c)$ of Proposition 1.4, any MASA in a U$_2$-factor is purely non-separable. 
\hfill $\square$

\vskip.05in

Let us also mention that   it was shown in (2.3.1$^\circ$ $(c)$ of [P13a]) 
that the Kadison-Singer paving problem over a MASA in a factor satisfying the stable U$_2$-property reduces to paving of projections 
having scalar expectation on the MASA.  
(Note that by Theorem 3.3 in [PV14], in order for a MASA $A$ in a II$_1$ factor $M$ to have the 
paving property, it is necessary that $A$ be purely non-separable.) Whether U$_2$ is a stable property was however left open in [P13a],     
but upon reading a preliminary draft of our paper  
Adrian Ioana pointed out to us that an argument  in the same vein as the proof of Theorem 1.2 easily implies U$_2$-stability as well.   
We thank him for sharing this with us.

\proclaim{1.6. Proposition} Condition $\text{\rm U}_2$ is a stable property.  
\endproclaim
\noindent
{\it Proof}. Assume the II$_1$ factor $M$ satisfies U$_2$. Since this 
 trivially  implies $\Bbb M_n(M)$ satisfies U$_2$, $\forall n$, to prove the stability   
it is sufficient to show that $pMp$ satisfies U$_2$ for any projection $p\in M$. 
Let $A_1, A_2\subset pMp$ be separable diffuse abelian von Neumann algebras. Let $R\subset M$ be a copy of the hyperfinite 
II$_1$ factor with $D\subset R$ its Cartan subalgebra and so that $p\in D$. Let also $\tilde{A}_i \subset M, i=1,2$, be separable diffuse abelian 
von Neumann algebras containing $p$ and such that $\tilde{A}_ip=A_i$. By the U$_2$ property of $M$, there exist unitaries $u_i\in M$ 
such that $u_i\tilde{A}_iu_i^*=D$. Since $D\subset M$ is Cartan, there exist $v_i\in \Cal N_R(D)$ such that 
$v_i(u_ipu_i^*)v_i^*=p$, $i=1,2$. But this means $w_i=v_iu_ip$ are unitaries in $pMp$ that conjugate $A_i$ onto $Dp$, $i=1,2$. 
Thus, $A_1, A_2$ are unitary conjugate as well. 
\hfill $\square$

\proclaim{1.7. Corollary} If a $\text{\rm II}_1$ factor $M$ satisfies property $\text{\rm U}_2$  
$($e.g.,  if $M$ is a $\text{\rm U}$-factor$)$, then a MASA $A\subset M$ has the paving property if and only if  any projection $q\in M$ with 
$E_A(q) \in \Bbb C1$ can be paved. 
\endproclaim
\noindent
{\it Proof.} By 1.6 above, property U$_2$ is stable, 
so the statement follows from (Proposition 2.3.1$^\circ$ (c) in [P13a]). 
\hfill $\square$

\vskip.05in
\noindent
{\bf 1.8. Remark.} While U$_1$ trivially implies U$_2$, we have no examples of a II$_1$ factor satisfying U$_2$ but not U$_1$. 
Note in this respect that if $M$ satisfies property U$_2$ and $A_0, A_1\simeq L^\infty[0,1]$ 
are von Neumann subalgebras of $M$ then by  conjugating by a unitary in $M$ 
we may assume $A_0=A_1$ and then property U$_1$ amounts to whether any automorphism of $(A_0, \tau)$ is implemented 
by a unitary in $M$. Thus, the following two additional properties of a II$_1$ factor $M$ are relevant: 

\vskip.05in
U$_3$  Given any separable diffuse abelian von Neumann subalgebra $A_0\subset M$, any automorphism of $(A_0, \tau)$ is implemented 
by a unitary in $M$.

\vskip.05in
U$'_3$  There exists a separable diffuse abelian von Neumann subalgebra $A_0\subset M$ such that any automorphism of $(A_0, \tau)$ is implemented 
by a unitary in $M$.  

\vskip.05in 

Thus, we see that U$_1 \Rightarrow \text{\rm U}_3\Rightarrow \text{\rm U}_3'$, U$_1 \Leftrightarrow (\text{\rm U}_2 + \text{\rm U}_3') \Leftrightarrow (\text{\rm U}_2 + \text{\rm U}_3)$, and that both U$_3$, U$'_3$ are stable properties (proof being similar to the proof of the stability of U$_1$, U$_2$). 
Thus, an example of a II$_1$ factor $M$ satisfying U$_2$ but not U$_1$  (so $M$ not a U-factor) should contain a copy of 
the non-atomic probability space $([0,1], \lambda)$ whose normalizer in $M$ does not implement all of its automorphism group.

\heading 2. Constructing disjoint singular MASAs in U-factors \endheading 

We show in this section that, under the continuum hypothesis, the size of the singular abelian core of any U-factor  is  quite ``large'' and can be estimated.  

We briefly recall (see e.g., [P81b]) that if $M$ is a II$_1$ factor and $A\subset M$ is a MASA, then $A$ is singular in $M$ iff any partial isometry $v\in M$ satisfying 
$v^*v, vv^*\in A$, $vAv^*\subset A$ must be contained in $A$. Also, using notations from intertwining theory (see e.g., 1.5 in [P13b], for 1.3 in [P16]) 
given two MASAs $A_1, A_2\subset M$ one has 
$A_1\prec_M A_2$ iff there exists a non-zero partial isometry $v\in M$ such that $v^*v\in A_1, vv^*\in A_2$ and $vA_1v^*\subset A_2$. 

\proclaim{2.1. Theorem} Let $M$ be a $\text{\rm II}_1$ $\text{\rm U}$-factor $M$ with the property that the cardinality 
of its unitary group $\Cal U(M)$ is equal to $\frak c$. If  the continuum hypothesis, 
$\frak c=\aleph_1$, is assumed, then $M$ contains more than $\frak c$ many mutually disjoint singular MASAs, i.e.,  
$\text{\rm r}(M)> \frak c$. Moreover, if the strong continuum hypothesis  $2^{\frak c}=\aleph_2$ is assumed, then ${\text{\rm r}}(M) = 2^{\frak c}$.
\endproclaim
\noindent 
{\it Proof}. Denote by $(I, <)$ the set of ordinals $<\aleph_1=\frak c$ endowed with its well ordered relation. Since $|\Cal U(M)|=\frak c$, 
it follows that $|\Cal P(M)|=\frak c$, and thus the cardinality of the set $\Cal V=\Cal V(M)=\{up \mid u\in \Cal U(M), p\in \Cal P(M)\}$ of partial isometries of $M$ 
is equal to $\frak c$ as well. Let $\{v_i\}_{i\in I}$ be an enumeration with repetition of $\Cal V$,  where each $v\in \Cal V$ appears $\frak c$-many times. 

Let $\Cal A$ be a maximal family of disjoint singular abelian wo-closed subalgebras $A \subset 1_AM1_A$ (which apriori may be an empty set). 
Assume $|\Cal A| \leq \frak c=\aleph_1$. Let $\{A_i\}_{i\in I}$ be a family of MASAs in $M$ indexed by our set $I$, 
such that each $A\in \Cal A$ appears as a direct summand of some $A_i$. 

Note that if we can show that under  these assumptions there exists a singular MASA $B\subset M$ such that 
$B\not\prec_M A_i$, $\forall i\in I$, then this would contradict the fact that $\{A_i\}_{i\in I}$ contains all of $\Cal A$, which  
was chosen to be the maximal singular core for $M$. This contradiction would show that one necessarily have $|\Cal A|>\frak c$, 
thus finishing the proof of the first part. If in addition we have $2^{\frak c}=\aleph_2$, since the total number of distinct MASAs in a II$_1$ factor $M$ 
with $|\Cal U(M)|=\frak c$ is obviously majorized by $2^{\frak c}$, it would then also follow that under the assumption $2^{\frak c}=\aleph_2$ 
one gets r$(M)=|\Cal A|=2^{\frak c}$. 

We construct $B$ as the wo-closure of the union of an increasing family $\{B_i\}_{i\in I}$ of separable diffuse abelian von Neumann subalgebras of $M$, which we 
construct by transfinite induction over $i \in I$, in the following way. 

Assume that $B_j$ have been constructed for all $j<i$. We want to construct $B_i$ 
so that $v_i$ is not intertwining $B_i$ into $B_i'\cap M$, nor $B_i$ into $A_j$ for $j\leq i$. To this end, we proceed as follows: 

$(a)$ Denote  $B^0_i=\overline{\cup_{j<i} B_j}$. Note that $B_i^0$ is separable abelian diffuse.  

$(b)$ If $v_i^*v_i \not\in B_i^0$ then by U$_0$ there exists a self-adjoint element $a\in (B_i^0)'\cap M$ such 
that $[v_i^*v_i, a]\neq 0$ and we let $B_i=B^0_i \vee \{a\}$. Note that $B_i$ is then still separable abelian and 
$[v_i^*v_i, B_i]\neq 0$. 

$(c)$ If $v_i^*v_i \in B^0_i$ then we let $K_i=\{j \in I_0, j\leq i \mid v_iB^0_iv_i^*\not\subset A_j\}$ and $L_i=\{j\in I_0, j\leq i \mid v_iB^0_iv_i^*\subset  A_j\}$. Note that 
$K_i, L_i$ are disjoint,  countable sets, with $K_i \cup L_i =\{j \in I_0 \mid  j \leq i\}$. Denote $p_i=v_i^*v_i\in B^0_i$ and notice that 
for each $j\in L_i$ we have $v_i^*A_jv_i \subset Q^0_i\overset{def}\to{=}(B^0_ip_i)'\cap p_iMp_i$, with $v_i^*A_jv_i$ a MASA in $Q^0_i$. Thus, if we denote 
$S_i:=\cup_{j\in L_i} v_i^*A_jv_i$ then the set  $S_i \subset Q^0_i$ is a countable union of abelian von Neumann algebras (even MASAs) in 
the II$_1$ von Neumann algebra $Q^0_i$, so $Q^0_i\setminus S_i$ is a $G_\delta$ dense subset of $Q^0_i$. 

Note already 
that if $a_0\in Q^0_i \setminus S_i$ is a self-adjoint element then any separable abelian von Neumann algebra that contains the abelian 
algebra $B^1_i=B^0_i \vee \{a_0\}$ cannot be intertwined by the partial isometry $v_i$ into $A_j$ for any $j\leq i$. 

In order to choose $B_i\supset B^1_i$ so that to exclude $v_i$ from properly normalizing any MASA $B$ containing $B_i$, let us note that there are 
several possibilities. 

$(i)$ $v_i\in B^1_i$, in which case we just put $B_i=B^1_i$. 

$(ii)$ $v_iB_i^1v_i^* \not\subset (B^1_i)'\cap M$, in which case we again let $B_i=B^1_i$. 

$(iii)$ $v_iB_i^1v_i^* \subset (B^1_i)'\cap M$ but $v_iB_i^1 v_i^* \not\subset B^1_i$. This means there exists $a\in B_i^1p_i$ such 
that $v_i a v_i^* \in ((B^1_i)'\cap M)\setminus B^1_i$, and by applying U$_0$  there  exists $a_1=a_1^*\in (B^1_i)'\cap M$ such that $[a_1, v_iav_i^*]\neq 0$. 
We then let $B_i=B^1_i \vee \{a_1\}$. 

$(iv)$ $v_iB_i^1 v_i^* \subset B^1_i$ but $v_iB_i^1 v_i^* \neq B^1_i v_iv_i^*$. In this case we have that $v_i^*B^1_iv_i$ strictly contains 
$B^1_i v_i^*v_i$. Like in $(iii)$ above, by U$_0$ there exist  $a'\in B^1_i$ and a self-adjoint $a_1'\in (B_i^1)'\cap M$ 
such that $[v_i^*a' v_i, a'_1]\neq 0$. We then define $B_i=B^1_i \vee \{a_1'\}$.

$(v)$ $v_iB_i^1 v_i^* = B^1_i v_iv_i^*$ but $v_i\not\in B^1_i$. This implies the partial isometry $v_i$ normalizes the II$_1$ von Neumann 
algebra $Q_i=(B^1_i)'\cap M$, acting non-trivially on it, having left and right supports in $\Cal Z(Q_i)=B^1_i$. There are two possibilities: 

$(v1)$ $v_i\in Q_i$. In this case $v_i^*v_i=v_iv_i^*=p_i \in \Cal Z(Q_i)$ and so $v_i$ is a non-central unitary in 
the II$_1$ von Neumann algebra $Q_ip_i$. 

We claim that if this is the case, then there exists a unitary $u\in Q_ip_i$ such that $v_iuv_i^*$ doesn't commute with $u$. 

To see this, first note that by Proposition 1.5 $(b)$, $Q_ip_i$ is of type II$_1$, so 
$Q_ip_i \not\prec_N \Cal Z(Q_ip_i)$  in any ambient II$_1$  factor $N$ that we would embed $Q_ip_i$. 
Taking $N$ to be a free product of $Q_ip_i$ with a diffuse tracial algebra, we can assume 
$Q_ip_i$ is embedded in a II$_1$ factor $N$ so that its relative commutant in $N$ is equal to $\Cal Z(Q_ip_i)$. But then 
 we can apply (Theorem 0.1 $(a)$ in [P13b]) to get a Haar unitary $u\in Q_ip_i$ that's approximately free to $x=v_i-E^N_{\Cal Z(Q_ip_i)}(v_i)\neq 0$. 
 In particular, one can take $u$ to be $\varepsilon$ $4$-independent to $x$, which for $\varepsilon>0$  sufficiently small insures that 
 $[v_iuv_i^*, u]\neq 0$. 
 
 Taking now $u\in Q_ip_i$ to be any unitary satisfying this property, we define $B_i=B^1_i \vee \{u\}$. 

$(v2)$ $v_i \not\in Q_i$. In this case $v_i$ acts non-trivially on the center of $Q_i$, so there exists  mutually orthogonal projections $z_1, z_2\in  \Cal Z(Q_i)$ 
such that $z_1 \leq v_i^*v_i$, $z_2\leq v_iv_i^*$ and $v_iz_1v_i^*=z_2$. Since $Q_iz_1$ is II$_1$, there exists a copy of $\Bbb M_2(\Bbb C)$ 
inside it. So there exist self-adjoint unitaries $u, w\in Q_iz_1$ such that $uw=-wu$. Let $c=u + v_i wv_i^*$ and define $B_i = B^1_i\vee \{c\}$. 
Note that $c, z_1, z_2$ are  elements in $B_i$ such that  $[v_i (cz_1)v_i^*, cz_2]\neq 0$. 

Finally, we define $B=\overline{\cup_i B_i}^{wo}$. Let us first show that $B$ is a MASA in $M$, i.e., $B=B'\cap M$. 
To see this, it is sufficient to prove that any selfadjoint unitary $v\in B'\cap  M$ lies in $B$. Since $v\in \Cal V$, it is of the form $v_i$ for some $i\in I$. 
This means $v_i$ is being considered in step $i$ of the induction and we see that we are necessarily in the situation $(v1)$, where we have chosen $B_i$ 
(which is a subalgebra of $B$) so that to contain some $b$ such that $v_ibv_i^*b\neq b v_ibv_i^*$, contradicting $[B, v_i]=0$. 

Assume now that $B$ is not singular. This implies there exists a non-zero partial isometry $w\in M$ with  $w^*w, ww^*$ mutually orthogonal 
projections in $B$. Thus $w\in \Cal V$ so $w=v_i$ for some $i\in I$ and so we have considered $w$ at step $i$ of the induction, and we are 
necessarily in one of the situations $(iii), (iv), (v1), (v2)$, which all lead to contradictions. 

Finally, assume $B\prec_M A_j$ for some countable ordinal $j\in I$. This means there exists a partial isometry $v\in M$ such that $v^*v\in B$, $vv^*\in A_j$ 
and $vBv^*=A_jvv^*$. Because of our choice of repeating $v$ $\frak c$-many times in $\{v_i\}_{i\in I}$, there exists $i\in I$ such that $i>j$ and $v=v_i$. 
But then the choices we made in $(b), (c)$ for the algebra $B_i\subset B$, easily imply that we cannot have $v_iBv_i^*\subset A_i$. 
\hfill $\square$

\proclaim{2.2. Corollary}  Let $\{M_n\}_{n\geq 1}$  be a sequence of separable tracial factors with $\text{\rm dim}(M_n) \rightarrow \infty$ and $\omega$ 
a free ultrafilter on $\Bbb N$. Denote   $M=\Pi_\omega M_n$ the associated ultraproduct $\text{\rm II}_1$ factor. If we assume the continuum 
hypothesis then 
$\text{\rm r}(M)> \frak c$. If we further assume the strong continuum hypothesis, then $\text{\rm r}(M)=2^{\frak c}$. 
\endproclaim
\noindent 
{\it Proof}. Since any ultraproduct II$_1$  factor $M=\Pi_\omega M_n$ satisfies the bicommutant axiom U$_0$, it is a U-factor. 
If in addition $M_n$ are all separable, then $|\Cal U(M_n)|=\frak c$, so $|\Cal U(M)|=\frak c^{\aleph_0}=\frak c$. Thus, 
we can apply Theorem 2.1 to conclude that under the CH condition we have r$(M)> \frak c$. Since the total number of distinct MASAs 
in $M$  is majorised  by the number of subsets of $\Cal U(M)$, it is bounded by $2^{\frak c}$. Thus, r$(M)\leq 2^{\frak c}$. 
So, if SCH is assumed then r$(M)=2^{\frak c}$.
\hfill $\square$

\heading 3. Further considerations  \endheading

The motivation behind our calculations of singular abelian rank of ultraproduct II$_1$ factors 
was  the hope that this invariant might be able to  differentiate among some of these factors (for instance, between    
$\Pi_\omega \Bbb M_{k_n}(\Bbb C)$, with  $k_n \nearrow \infty$, and $M^\omega$, for a separable non-Gamma II$_1$ factor $M$). 
But our calculations, which anyway depend on CH/SCH, show that, like  in the separable 
case where one has r$(M)=\frak c$ for any separable II$_1$ factor  $M$ (cf. [P13a], [P18]; see Remark 2.7 in [BDIP23]), 
the singular abelian rank is the same, equal to $2^{\frak c}$, for all ultraproducts II$_1$ factors. 

One can try to ``diminish'' the number of disjoint singular MASAs by restricting our attention to MASAs that satisfy various 
stronger versions of singularity, thus attempting to bring them to a ``small cardinality'', even finite if possible. 

Thus, in the spirit of the terminology in (Definitions 2.5, 2.9 in [BDIP23]), let us denote by $\Cal A^*_M$ 
a maximal family of disjoint ``special'' singular MASAs in the II$_1$ factor $M$ satisfying a ``generic'' stronger 
singularity property $^*$. As in ([BDIP23]), we will in fact view $\Cal A^*_M$  in ``unfolded'' form, 
as one single singular abelian wo-closed $^*$-subalgebra  generated by finite projections 
in the II$_\infty$ factor $\Cal M=M\overline{\otimes} \Cal B(\ell^2K)$, where $K$ is a set of  sufficiently large cardinality ($K\geq 2^{|\Cal U(M)|}$ will do),  
which is so that any of its finite corners has the property $^*$, and which is maximal 
(with respect to inclusion) with these properties. Note that these requirements force the definition of disjointness  to be taken possibly stronger as well.  

One  then takes the corresponding rank r$_*(M)$ to be the trace $Tr_\Cal M$ of the support of 
$\Cal A^*_M\subset \Cal M$. Like in ([BDIP23]) one clearly has the amplification formula r$_*(M^t)=\text{\rm r}_*(M)/t$, $\forall t>0$, making such  
considerations particularly interesting if the rank of the ``special'' singular core could be shown finite. 

We illustrate below with four examples of  such a possible strengthening.

\vskip.05in
\noindent
{\it $3.1$. The supersingular abelian core}. Following ([P13a]),  we'll say that a wo-closed abelian 
$^*$-subalgebra $A$ in a II$_1$ factor $M$ is {\it supersingular} if there is no automorphism $\theta\in \text{\rm Aut}(M)$ 
such that $\theta(Ap)\subset A$ for some non-zero $p\in \Cal P(A)$ other than the inner automorphisms of $M$ that 
act trivially on $pMp$. Two such supersingular abelian subalgebras $A_1, A_2 \subset M$ are {\it disjoint} 
if there exists no automorphism $\theta$ of $M$ satisfying $\theta(A_1p_1)\subset A_2$ for some non-zero projection 
$p_1\in A_1$. Note that this is the same as requiring that $A_1 \oplus A_2$ be supersingular in $M^2=\Bbb M_2(M)$. 

As we mentioned above, like in ([BDIP23]), we in fact view any family $\Cal A$ of disjoint (in this stronger sense) supersingular abelian subalgebras in $M$ 
in its ``unfolded'' form, as one single supersingular abelian algebra generated by finite projections in $M\overline{\otimes} \Cal B(\ell^2K)$, 
for a sufficiently large $K$.  One clearly has a maximal  such algebra with respect to inclusion,  $\Cal A_M^{\text{\rm ss}}$, 
which is moreover unique up to   
unitary conjugacy in $\Cal M$, and which we'll call the {\it superrsingular abelian core}. The corresponding  {\it supersingular rank} r$_{\text{\rm ss}}(M)$ 
is then given by the trace $Tr_\Cal M$ of the support of $\Cal A_M^{\text{\rm ss}}$  in $\Cal M$, viewed as a cardinality when infinite.

\vskip.05in
\noindent
{\it $3.2$. The coarse abelian core}. In the same spirit, this time following ([P18]), 
one can take in $\Cal M=M\overline{\otimes}\Cal B(\ell^2K)$ the {\it coarse abelian core} 
to be a wo-closed abelian $^*$-subalgebra $\Cal A^{\text{\rm c}}_M\subset \Cal M$ generated by finite projections 
with the property that $\Cal Ap$ is coarse in $p\Cal Mp$ for any finite projection 
$p\in \Cal A^{\text{\rm c}}_M$, and which is maximal with respect to inclusion. Note that disjointness for coarse 
abelian $A_1, A_2\subset M$ amounts to $A_1, A_2$ being a coarse pair (as defined in [P18]).

The coarse core this way defined is clearly unique in $\Cal M$ up to 
unitary conjugacy. The {\it coarse abelian rank} is then 
$\text{\rm r}_{\text{\rm c}}(M)=Tr_{\Cal M}(1_{\Cal A^{\text{\rm c}}_M})$.  

Note however that by results in [P18], for any separable $M$ one has $\text{\rm r}_{\text{\rm c}}(M)>\aleph_0$, 
so if we assume CH then $\text{\rm r}_{\text{\rm c}}(M)=\frak c=\aleph_1$.

\vskip.05in
\noindent
{\it $3.3$. The maximal amenable abelian core.}  We define the {\it maximal amenable abelian core} $\Cal A_M^{\text{\rm ma}}$ of 
the II$_1$ factor $M$ as the  $\Cal M=M\overline{\otimes}\Cal B(\ell^2K)$ 
the wo-closed abelian $^*$-subalgebra $\Cal A=\Cal A^{\text{\rm ma}}_M\subset \Cal M$ generated by finite projections 
with the property that $\Cal A$ is maximal amenable in $\Cal M$, and which is maximal with respect to inclusion. 
Its {\it maximal amenable abelian rank} by $\text{\rm r}_{\text{\rm ma}}(M)=Tr_{\Cal M}(1_{\Cal A^{\text{\rm ma}}_M})$. 

While it is not clear how this invariant fares for separable II$_1$ factors, note that by (Theorem 5.3.1 in [P13a]) any ultraproduct  $A=\Pi_\omega A_n$ of singular MASAs  in II$_1$ factors $A_n\subset M_n$,  is maximal amenable in $M=  \Pi_\omega M_n$. Thus, for such factors one has r$_{\text{\rm ma}}(M)\geq \frak c$.

\vskip.05in
\noindent
{\it $3.4$. The singular s-MASA core.}  Following ([P16]), a MASA  $A$ in a II$_1$ factor  $M$ is an $s$-MASA if 
$A\vee A^{op}$ is a MASA in $\Cal B(L^2M)$. By a well know result of Feldman and Moore ([FM77]), 
any Cartan subalgebra satisfies this property. It has been shown in ([P16]) that  if the II$_1$ factor $M$ is separable and has s-MASAs, then it has 
singular s-MASAs, and in fact it has $> \aleph_0$ many disjoint s-MASAs. 

One defines the {\it s-MASA core} of a II$_1$ factor $M$, as the wo-closed abelian $^*$-subalgebra $\Cal A=\Cal A^{\text{\rm s}}_M\subset \Cal M$ generated by finite projections 
with the property that $\Cal Ap$ is a singular s-MASA in $p\Cal Mp$ for any finite projection 
$p\in \Cal A$, and which is maximal with respect to inclusion. Again, this is obviously unique in $\Cal M$ up to 
unitary conjugacy. The {\it s-MASA rank} of $M$ is then 
$\text{\rm r}_{\text{\rm s}}(M)=Tr_{\Cal M}(1_{\Cal A^{\text{\rm s}}_M})$. So by ([P16]), in this case as well the associated rank is huge, r$_{\text{\rm s}}(M) > \aleph_0$, 
so equal to $\frak c$ when CH is assumed. It is not clear if ultraproduct factors, or even more generally U-factors, can have singular s-MASAs at all. 
Since existence of  an  s-MASA in a II$_1$ factor is a ``thinness'' property that ultraproducts are unlikely to have,  
it seems that such factors cannot have s-MASAs, but this remains an open problem.

\vskip.05in
\noindent
%{\it $3.5$. Freely complemented MASAs.}  

\head  References \endhead

\item{[AP17]} C. Anantharaman, S. Popa: ``An introduction to II$_1$ factors'', \newline www.math.ucla.edu/$\sim$popa/Books/IIun-v13.pdf 

\item{[BDIP23]} R. Boutonnet, D. Drimbe, A. Ioana, S. Popa: {\it Non-isomorphism of $A^{*n}, 2\leq n \leq \infty$, for a non-separable abelian von Neumann algebra $A$}, 
Geometric and Functional Analysis (GAFA), {\bf 34} (2024), arXiv:2308.05671.

\item{[D54]} J. Dixmier: {\it Sous-anneaux ab\'eliens maximaux dans les facteurs de type fini}, Ann. of Math. {\bf 59} (1954), 279-286. 

\item{[D57]} J. Dixmier: ``Les alg\'ebres d'op\'erateurs sur l'espace Hilbertien (Alg\'ebres de von Neumann)'', Gauthier-Villars, Paris, 1957.

\item{[D64]} J. Dixmier, C. Lance: ``Les C$^*$-algebres et leurs representations'', Gauthier-Villars, Paris, 1964. 

\item{[FM77]} J. Feldman, C. Moore: {\it Ergodic equivalence relations, cohomology, and von Neumann algebras}, 
II. Trans. Amer. Math. Soc., {\bf 234(} (1977), 325Ð359.

\item{[GKPT24]} D. Gao, S. Kunnawalkam Elayavalli, G. Patchell, H. Tan: {\it A highly indecomposable} II$_1$ {\it factor}, in preparation. 

\item{[HP17]}  C. Houdayer, S. Popa:  {\it Singular MASAs in type} III {\it factors and Connes' bicentralizer problem}, 
Proceedings of the 9th MSJ-SI "Operator Algebras and Mathematical Physics" held in Sendai, Japan, 2016 (math.OA/1704.07255)   

%\item{[I24]} A. Ioana: Private communication, February 2024. 

\item{[P81a]} S. Popa: {\it Orthogonal pairs of *-subalgebras in
finite von Neumann algebras}, J. Operator Theory, {\bf 9} (1983),
253-268.

\item{[P81b]} S. Popa: {\it Singular maximal abelian $^*$-subalgebras in continuous von Neumann algebras}, Journal of Funct. Analysis, 
{\bf 50} (1983), 151-166. 

\item{[P13a]} S. Popa: {\it A} II$_1$ {\it factor approach to the Kadison-Singer problem}, 
Comm. Math. Physics. {\bf 332} (2014), 379-414 (math.OA/1303.1424).

\item{[P13b]} S. Popa: {\it Independence properties in subalgebras of ultraproduct} II$_1$ {\it factors}, Journal of Functional Analysis 
{\bf 266} (2014), 5818Ð5846 (math.OA/1308.3982). 

\item{[P16]}  S. Popa:  {\it Constructing MASAs with prescribed properties}, Kyoto J. of Math,  {\bf 59} (2019), 367-397 (math.OA/1610.08945). 

\item{[P18]} S. Popa: {\it Coarse decomposition of} II$_1$ {\it factors}, Duke Math. J. {\bf 170} (2021) 3073 - 3110 (math.OA/1811.09213)

\item{[P20]} S. Popa: ``Topics in II$_1$ factors'', graduate courses at UCLA, Winter 2020, Fall 2022. 

\item{[PV14]}  S. Popa, S. Vaes: {\it Paving over arbitary MASAs in von Neumann algebras}, 
Analysis and PDE  (2015) 101-123 (math.OA/1412.0631)

\item{[Pu60]} L. Pukanszky: {\it On maximal abelian subrings in factors of type} II$_1$, Canad. J. Math. {\bf 12} (1960), 289-296.

\item{[R91]} F. Radulescu: {\it Singularity of the radial subalgebra of $L(\Bbb F_N)$ and the Pukanszky invariant}, Pacific J. Math., {\bf 151} (1991), 297-306.
 
\enddocument